\documentclass[namedate,webpdf,imanum]{ima-authoring-template}

\usepackage{bm}
\usepackage{subfig}
\usepackage{hyperref}
\usepackage{MnSymbol}
\usepackage{booktabs}
\usepackage{caption}
\usepackage{amsmath,amsthm,amsfonts}
\usepackage{etoolbox}
\usepackage{cprotect}
\usepackage{cleveref}
\usepackage{mathtools}

\crefname{algorithm}{Algorithm}{Algorithms}
\crefname{chapter}{Chapter}{Chapters}
\crefname{equation}{}{Equations}
\crefname{figure}{Fig.}{Figures}
\crefname{table}{Table}{Tables}
\crefname{lemma}{Lemma}{Lemmas}
\crefname{section}{Section}{Sections}
\crefname{subsection}{Section}{Sections}

\graphicspath{{figures/}}

\theoremstyle{thmstyletwo}%
\newtheorem{theorem}{Theorem}
\newtheorem{remark}{Remark}%

\newtheorem{assumption}{Assumption}

\numberwithin{equation}{section}

\def\bu{\boldsymbol{u}}
\def\md{\mathrm{d}}
\def\mB{\mathcal{B}}
\def\mD{\mathcal{D}}
\def\mF{\mathcal{F}}
\def\mG{\mathcal{G}}
\def\mL{\mathcal{L}}
\def\mP{\mathcal{P}}
\def\mS{\mathcal{S}}
\def\mT{\mathcal{T}}
\def\mU{\mathcal{U}}
\def\mX{\mathcal{X}}

\begin{document}

\DOI{}
\copyrightyear{}
\vol{}
\pubyear{}
\access{Advance Access Publication Date: Day Month Year}
\appnotes{Paper}
\copyrightstatement{Published by Oxford University Press on behalf of the Institute of Mathematics and its Applications. All rights reserved.}
\firstpage{1}


\title[Optimal Spectral Poisson Solver]{An optimal complexity spectral solver for the Poisson equation}

\author{Ouyuan Qin*
\address{\orgdiv{School of Mathematical Sciences}, \orgname{University of Science and Technology of China}, \orgaddress{\street{96 Jinzhai Road}, \state{Hefei}, \postcode{Anhui 230026}, \country{China}}}}
\authormark{O.Qin}

\corresp[*]{Corresponding author: \href{email:ouyuanqin@mail.ustc.edu.cn}{ouyuanqin@mail.ustc.edu.cn}}

\received{Date}{0}{Year}
\revised{Date}{0}{Year}
\accepted{Date}{0}{Year}


\abstract{We propose a spectral solver for the Poisson equation on a square domain, achieving optimal complexity through the ultraspherical spectral method and the alternating direction implicit (ADI) method. Compared with the state-of-the-art spectral solver for the Poisson equation \citep{for}, our method not only eliminates the need for conversions between Chebyshev and Legendre bases but also is applicable to more general boundary conditions while maintaining spectral accuracy. We prove that, for solutions with sufficient smoothness, a fixed number of ADI iterations suffices to meet a specified tolerance, yielding an optimal complexity of $\mathcal{O}(n^2)$. The solver can also be extended to other equations as long as they can be split into two one-dimensional operators with nearly real and disjoint spectra. Numerical experiments demonstrate that our algorithm can resolve solutions with millions of unknowns in under a minute, with significant speedups when leveraging low-rank approximations.}


\maketitle

\section{Introduction}\label{sec:intro}
The Poisson equation is a simple but fundamental problem in partial differential equations (PDEs) and frequently used as a benchmark for various numerical methods. Consider the Poisson equation with homogeneous Dirichlet boundary conditions on a square domain:
\begin{align}
  \label{poissonD}
  u_{xx} + u_{yy} = f, \quad (x, y) \in [-1,1]^2, \quad u(\pm 1, \cdot) = u(\cdot, \pm 1) = 0,
\end{align}
where $f$ is a known continuous function and $u$ is the solution we seek. A general approach to solving \eqref{poissonD} involves constructing differential matrices $D_{xx}$ and $D_{yy}$ for the $x$- and $y$-directions, and then solving the Sylvester matrix equation
\begin{align}
  \label{fdsylvester}
  X D_{xx}^T + D_{yy} X = F,
\end{align}
where $X , F \in \mathbb{R}^{n \times n}$ are the corresponding discretizations of $u$ and $f$, respectively. The classical method for solving \eqref{fdsylvester}, i.e., Sylvester equation, is the Bartels--Stewart (B--S) algorithm \citep{bar}, which requires Schur or eigenvalue decompositions of $D_{xx}$ and $D_{yy}$.

For the finite difference (FD) method, the differential matrices in \eqref{fdsylvester} are constructed as tridiagonal Toeplitz matrices for a 5-point stencil on an equispaced grid. The eigenvectors of these matrices are exactly the discrete Fourier transform matrix \citep{gol2, lev}. With the aid of fast Fourier transform (FFT) \citep{lev}, equation \eqref{fdsylvester} can be solved in $\mathcal{O}(n^2 \log n)$ operations. However, the discretization error of FD method is relatively large, and high-accuracy solutions are better achieved using spectral methods.

When spectral method is considered, the eigendecomposition of its differentiation matrices does not exhibit obvious structure and Schur decomposition is required for solving \eqref{fdsylvester} \citep{hai, she2}, resulting an $\mathcal{O}(n^3)$ complexity algorithm. Recently, a fast Poisson spectral solver has been proposed \citep{for}, which takes advantage of the spectra of differential matrices instead of full eigendecompositions. This method reports an asymptotic complexity of $\mathcal{O}(n^2 (\log n)^2)$ and is applicable to various domains. The solver employs carefully designed bases for \eqref{poissonD} and results in symmetric pentadiagonal coefficient matrices. Despite the involvement of the alternating direction implicit (ADI) method for solving \eqref{fdsylvester}, the algorithm in \citet{for} can be seen as a direct method, with guarantees provided by rational approximations and bounds on the Zolotarev number \citep{leb, lu}. Numerical results demonstrate that this method is faster than any available direct methods. This approach has also been extended to the $hp$-finite element method ($hp$-FEM) \citep{kno}, where the coefficient matrices are banded-block-banded-arrowhead. This method can handle discontinuities in the input data and achieve both spectral convergence and quasi-optimal complexity.

Unfortunately, the fast solver in \citet{for} is limited to Poisson equations with Dirichlet boundary conditions due to the special basis used. For Neumann and Robin boundary conditions, the problems must be solved using the ultraspherical spectral (US) method \citep{olv}. However, employing the techniques from US method causes the banded structure to be lost, and the spectrum of the US differential operator remains unclear. With the goal of proposing a unified Poisson spectral solver with optimal complexity, we introduce a new approach for \eqref{poissonD}, where the differentiation operators are replaced by operators from the US method, and the boundary conditions are satisfied through basis recombination \citep{qin2}. The resulting Sylvester equation is solved using the ADI method. A detailed analysis of the spectrum of the differential operator, as approximated by the US method with various boundary conditions, is presented here for the first time. This analysis addresses the unresolved questions raised in \citet[\S 6.2]{for} and guarantees a finite number of ADI iterations. Furthermore, while the number of iterations required to achieve a given tolerance grows as $\mathcal{O}(\log n)$ in \citet{for}, we show that $\mathcal{O}(1)$ iterations suffice, provided that the underlying solution is sufficiently smooth. Since each iteration of the ADI method requires $\mathcal{O}(n^2)$ operations, this result implies a truly optimal $\mathcal{O}(n^2)$ algorithm for solving \eqref{poissonD}. This promising result can be extended to other strongly elliptic PDEs, as long as the spectra of the operators in the $x$- and $y$-directions are nearly real and disjoint. Additionally, we provide several useful implementation details not mentioned in \citet{for}, such as error estimation during the solving process, recycling solutions for adaptivity, and optimizing the order of shifts. Finally, we present a low-rank version of our solver in case a low-rank factorization of the right-hand side $F$ is available. These seemingly minor aspects contribute to the efficiency and robustness of our solver.

The paper is organized as follows: \cref{sec:us} presents a brief review of the US method and the basis recombination technique. \cref{sec:adi} discusses the ADI method for solving the Sylvester equation, along with an estimation of the spectrum of the US differential operator and a crucial theorem that guarantees the sufficiency of fixed ADI iterations for solving PDEs. \cref{sec:details} provides additional implementation details that are highly beneficial for applying the new method. In \cref{sec:more}, we address general boundary conditions, separable coefficients, fourth-order equations, and nonhomogeneous boundary conditions. Numerical experiments in \cref{sec:exp} validate the proposed method in terms of both accuracy and speed, before we conclude in the final section.

Throughout this article, calligraphic fonts are used for operators or infinite matrices, and capital letters are used for finite matrices and truncations of operators. Infinite vectors and functions are denoted by bold and normal fonts, respectively. The standard 2-norm and its induced matrix norm are used, unless otherwise specified.

\section{Ultraspherical spectral method}\label{sec:us}
We start with one-dimensional problem for the very essence of US method. For $u_{xx} = f$ on $[-1, 1]$ with $u(\pm 1)=0$, the US method seeks a Chebyshev polynomial series $u(x) = \sum_{i=0}^{\infty} u_i T_i(x)$, where $T_i(x)$ is the Chebyshev polynomial of degree $i$. A sparse $\mD_{\ell}$ is utilized to represent $\md^{\ell}/\md x^{\ell}$, where
\begin{align}
  \label{usdiff}
  \begin{aligned}
    & ~~~\, \mathop{\rotatebox[origin=c]{180}{$\underbrace{\hphantom{00 \cdots 0}}$}}^{\ell \text{ times}} \\[-9pt]
    \mD_{\ell} = 2^{\ell - 1}(\ell - 1)!  & \begin{pmatrix}
        0 ~~ \cdots ~~ 0 &\ell & & & &  \\
        & &\ell+1& & & \\
        & & &\ell+2 & & \\
        & & & &\ddots& \\
    \end{pmatrix}
  \end{aligned}, \quad \ell \geq 1
\end{align}
maps the Chebyshev coefficients to the ultraspherical $C^{(\ell)}$ coefficients\footnote{Ultraspherical polynomials $C^{(\ell)}_i(x)$, $i=0,1,\dots$ are orthogonal on $[-1, 1]$ with respect to the weight function $(1-x^2)^{\ell-1/2}$ \citep{olv2}.}. If $f$ is also expressed as a Chebyshev series $\sum_{i=0}^{\infty} f_i T_i(x)$, conversion operators $\mS_1$ and $\mS_0$ should be applied to match the image space of $\mD_2$, where
\begin{align}
  \label{usconv}
  \mS_0 = \begin{pmatrix}
    1&  &-\frac{1}{2}& & \\
     &\frac{1}{2}& &-\frac{1}{2}& \\
     & &\frac{1}{2}& &\ddots \\
     & & &\ddots& \\
    \end{pmatrix} \text{ and }
  \mS_{\ell} = \begin{pmatrix}
    1&  &-\frac{\ell}{\ell +2}& & \\
     &\frac{\ell}{\ell +1}& &-\frac{\ell}{\ell +3}& \\
     & &\frac{\ell}{\ell +2}& &\ddots\\
     & & &\ddots& \\
    \end{pmatrix}, \quad \ell \geq 1
\end{align}
map the coefficients in Chebyshev and $C^{(\ell)}$ to those in $C^{(1)}$ and $C^{(\ell+1)}$, respectively. Thus, the equation $u_{xx} = f$ is equivalent to
\begin{align}
  \label{1dpoisson}
  \mD_2 \bu = \mS_1 \mS_0 \boldsymbol{f},
\end{align}
where $\bu = (u_0, u_1, \dots)^{\top}$ and $\boldsymbol{f} = (f_0, f_1, \dots)^{\top}$ are Chebyshev coefficients of $u$ and $f$, respectively.

To incorporate boundary conditions into \eqref{1dpoisson}, the boundary bordering technique is applied in \citet{olv}, where several dense rows representing boundary conditions are placed at the top of the linear system. In \citet{tow}, a reduced system of \eqref{1dpoisson} is constructed by removing degrees of freedom through boundary conditions. Both approaches for treating boundary conditions result in an almost-banded system. To preserve the banded structure of \eqref{usdiff} and \eqref{usconv}, we adopt another method, namely basis recombination. This approach has been proven successful in the literature \citep{she2,jul,qin2}, and many advantages of the US method can be retained with it.

For any linear constraints $\mB u = (\mB_1, \dots, \mB_N)^T u = (0, \dots, 0)^T$, we have
\begin{align*}
  \begin{pmatrix}
    0\\
    \vdots\\
    0
  \end{pmatrix} =
  \begin{pmatrix}
    \mB_1 \\
    \vdots \\
    \mB_N
  \end{pmatrix} u = 
  \begin{pmatrix}
    \mB_1 \\
    \vdots \\
    \mB_N
  \end{pmatrix} \left(\sum_{i=0}^{\infty} u_i T_i(x)\right) = 
  \begin{pmatrix}
    b_{10} & b_{11} & \cdots\\
    \vdots & \vdots & \vdots \\
    b_{N0} & b_{N1} & \cdots
  \end{pmatrix} \bu,
\end{align*}
where $b_{ji} = \mB_j(T_i(x))$ for $i=0,1,\dots$ and $j=1,\dots,N$. Denote the $i$-th column of the matrix above by $B_i = (b_{1i}, \dots, b_{Ni})^T$. For any recombined basis $v^k(x) = \sum_{i \in S_k} v_i^k T_i(x)$ satisfying boundary conditions $\mB v^k = (0, 0, \dots, 0)^T$, we know that $\sum_{i \in S_k} B_i v_i^k = (0, 0, \dots, 0)^T$, where $S_k \subset \{0,1,2,\dots\}$ is an index set. The choice of $S_k$ has a direct relation to the structure of the final linear system. Since we want to keep the system in \eqref{1dpoisson} banded after equipping it with boundary conditions, $S_k$ with aggregated elements is preferred \citep{she2, aur, qin2}, i.e., only several adjacent indices are contained in $S_k$. Considering that there are $N$ linearly independent constraints, $|S_k| > N$ is a necessary condition for the existence of nontrivial $v^k(x)$. Thus, we take $S_k = \{k, k+1, \dots, k+N\}$ and the system
\begin{align}
  \label{reducedBC}
  \begin{pmatrix}
    \mB_{k1} & \mB_{k+1,1} & \cdots & \mB_{k+N,1} \\
    \mB_{k2} & \mB_{k+1,2} & \cdots & \mB_{k+N,2} \\
    \vdots & \vdots & \ddots & \vdots \\
    \mB_{kN} & \mB_{k+1,N} & \cdots & \mB_{k+N,N} \\
  \end{pmatrix}
  \begin{pmatrix}
    v_{k}^k \\
    v_{k+1}^k \\
    \vdots \\
    v_{k+N}^k
  \end{pmatrix} = 
  \begin{pmatrix}
    0 \\
    0 \\
    \vdots\\
    0
  \end{pmatrix}
\end{align}
would admit a nonzero solution if one of $v_k^k, \dots, v_{k+N}^k$ is specified with a given value.

For example, if $\mB$ are Dirichlet conditions, i.e., $\mB_{1i} = (-1)^{i}$ and $\mB_{2i} = 1$, a particular solution to \eqref{reducedBC} is $v^k_k = 1$, $v^k_{k+1} = 0$, and $v^k_{k+2} = -1$, which is exactly the same basis $\phi^k(x) = T_k(x) - T_{k+2}(x)$ found in \citet{she2}. If we assemble coefficients of recombined bases column by column, a transformation operator $\mT_{D}$ related to homogeneous Dirichlet conditions is formed, where
\begin{align}
  \label{transD}
  \mT_{D} = \begin{pmatrix}
    \ \ 1 & & & \\
     & \ \ 1 & & \\
    -1 &  & \ \ 1 & \\
     & -1 & & \ddots \\
     &  & \ddots & \\
  \end{pmatrix}
\end{align}
maps the coefficients in recombined Chebyshev basis to those in canonical Chebyshev basis. Note that $\mT_{D}$ is lower triangular with a lower bandwidth of 2. With the transformation operator, the equation in \eqref{1dpoisson} together with homogeneous Dirichlet conditions becomes
\begin{align}
  \label{1dpoissonbc}
  \mD_2 \mT_{D} \tilde{\bu} = \mS_1 \mS_0 \boldsymbol{f},
\end{align}
where $\tilde{\bu}$ contains coefficients of the solution in the recombined Chebyshev basis, and it can be converted back to Chebyshev coefficients by $\bu = \mT_{D} \tilde{\bu}$.

With the ingredients for the one-dimensional problem in hand, the extension to the two-dimensional Poisson equation is an easy task using tensor products. Under the assumption that the solution of \eqref{poissonD} is a tensor product Chebyshev series, i.e.,
\begin{align}
  \label{usapprox}
  u(x, y) = \sum_{i=0}^{\infty}\sum_{j=0}^{\infty} u_{ij} T_i(y) T_j(x), \quad (x, y) \in [-1, 1]^2,
\end{align}
the $x$- and $y$-differentiations in \eqref{poissonD} become $\mU \mD_2^{\top}$ and $\mD_2 \mU$, where $\mU=(u_{ij})_{i,j=0}^{\infty}$ is an infinite matrix representing the coefficients of $u(x, y)$ in the tensor product Chebyshev basis. However, they cannot be added straightforwardly as their bases for the $x$- and $y$-directions are different. To get an equation in the tensor product ultraspherical polynomial basis, we apply conversion operators and basis recombination, arriving at
\begin{align}
  \label{poissonsylbc}
   \mS_1 \mS_0 \mT_{D} \tilde{\mU} (\mD_2 \mT_{D})^{\top} + \mD_2 \mT_{D} \tilde{\mU} (\mS_1 \mS_0 \mT_{D})^{\top} = \mS_1 \mS_0 \mF (\mS_1 \mS_0)^{\top},
\end{align}
where $\tilde{\mU}$ is the unknown in the tensor product of recombined basis and $\mF$ are the coefficients of the function $f(x, y)$ in the tensor product Chebyshev basis. Once $\tilde{\mU}$ is obtained, the usual tensor product Chebyshev coefficients are given by $\mU = \mT_{D} \tilde{\mU} \mT_{D}^{\top}$.

Although there are evident advantages of boundary bordering \citep{olv}, we choose basis recombination for several convincing reasons:
\begin{enumerate}
  \item A basis recombination with aggregated elements can always be sought, which leads to a banded structure in \eqref{poissonsylbc} ($\mD_2 \mT_{D}$ and $\mS_1 \mS_0 \mT_{D}$ are both banded) instead of the almost-banded structure generated in \citet{olv} and \citet{tow}. This allows for the use of highly optimized LAPACK \citep{and} subroutines to solve banded linear systems for one-dimensional problems.
  \item Despite the solution of \eqref{poissonsylbc} being a tensor product of recombined basis, it can be easily transferred back to the Chebyshev basis with the left and right application of triangular banded transformation operators, and the cost is negligible.
  \item The extra work compared to boundary bordering involves solving \eqref{reducedBC} and obtaining the formulas for the transformation operator $\mT_D$. The complexity of this step is $\mathcal{O}(N^3 n)$ if solved numerically or $\mathcal{O}(1)$ if solved analytically, which does not threaten the optimal complexity of the US method. A systemic way for determining the recombined basis is shown in \citet{qin2}.
  \item Note that the solutions of \eqref{reducedBC} form a one-dimensional subspace, and there is flexibility in choosing elements of $\mT_D$. Besides the simple one in \eqref{1dpoissonbc}, we could scale each column in $\mT_D$ so that the main diagonal of $\mD_{\ell} \mT$ are all 1s, making our method well-conditioned without the extra preconditioner used in \citet[\S 4]{olv}.
\end{enumerate}

\section{ADI method}\label{sec:adi}

Traditional methods for solving the Sylvester matrix equation, including \citet{bar, gol}, require costly Schur decompositions, and the accuracy of these methods depends heavily on the quality of the computed orthogonal matrices. While ADI has long been used as an iterative method for solving Sylvester equations \citep{bir, pea, wac}, it can be considered a direct method if the properties of separate spectra are exploited \citep{for} and no eigenvectors are explicitly involved in the iterations. Since our method also enjoys the properties introduced in \citet{for}, we employ ADI as the fast solver for \eqref{poissonsylbc} and elaborate on it below.

For simplicity, denote $D_{yy}$ and $D_{xx}^T$ in \eqref{fdsylvester} by $A$ and $B$ respectively, i.e.,
\begin{align}
  \label{modelsyl}
  AX + XB = F, \quad A, B, F \in \mathbb{R}^{n \times n},
\end{align}
where $X \in \mathbb{R}^{n \times n}$ is unknown. The ADI method generates new iterations by solving systems related to each direction sequentially. With shifts $p_j$ and $q_j$ selected for the $j$-th iteration, the approximation $X_j$ is computed as
\begin{align}
  \label{adiiter}
  \begin{aligned}
    (A - p_j I_n) X_{j-1/2} &= F - X_{j-1} (B + p_j I_n),\\
    X_{j}(B + q_j I_n) &= F - (A - q_j I_n) X_{j-1/2},
  \end{aligned}
\end{align}
where $X_0=0$ is the default initial iteration and $I_n$ is the square identity matrix of size $n$. Combining \eqref{modelsyl} and \eqref{adiiter}, it can be shown that the error satisfies \citep{sab}
\begin{align}
  \label{adierror}
  \begin{aligned}
    X - X_k &= (A - q_k I_n)(A - p_k I_n)^{-1} (X - X_{k-1}) (B + p_k I_n)(B + q_k I_n)^{-1} \\
            &= s_k(A) (X - X_0) s_k(-B)^{-1} = s_k(A) X s_k(-B)^{-1},
  \end{aligned}
\end{align}
where
\begin{align}
  \label{adisk}
  s_k(x) := \prod_{j=1}^{k} \frac{x - q_j}{x - p_j}.
\end{align}
If $A$ and $B$ can be diagonalized by $A = U_A \Lambda_A U_A^{-1}$ and $B = U_B \Lambda_B U_B^{-1}$, then
\begin{align}
  \label{relerror}
  \frac{\left\| X - X_k \right\|}{\left\| X \right\|} \leq \left\| s_k(A) \right\| \| s_k(-B)^{-1} \| \leq \kappa(U_A) \kappa(U_B) \max_{\substack{\lambda \in \sigma(A) \\ \mu \in \sigma(-B)}} \frac{|s_k(\lambda)|}{|s_k(\mu)|},
\end{align}
where $\kappa(V) = \| V \| \| V^{-1} \|$ and $\sigma(\cdot)$ denotes the spectrum of a matrix. An immediate observation is that $X_n = X$ if $p_j$ and $q_j$ are chosen as the eigenvalues of $-B$ and $A$, respectively. However, the exact spectrum of a matrix is difficult to obtain, and $n$ iterations of \eqref{adiiter} are too many to solve \eqref{modelsyl} speedily.

To achieve a minimum relative error, we should choose shifts $p_j$ and $q_j$ that minimize 
the last factor of the right-hand side in \eqref{relerror}. If the spectra of $A$ and $-B$ are real and disjoint, namely $\sigma(A) \subset [a, b]$ and $\sigma(-B) \subset [c, d]$ with $[a, b] \cap [c, d] = \varnothing$, the minimization can be relaxed by replacing the spectrum with the intervals as
\begin{align}
  \label{znumber}
  Z_k([a, b], [c, d]) := \min_{\substack{p_j \in \mathbb{C}\\ q_j \in \mathbb{C}}} \max_{\substack{\lambda \in [a, b]\\ \mu \in [c, d]}} \frac{|s_k(\lambda)|}{|s_k(\mu)|},
\end{align}
which is known as the Zolotarev number. This number has been extensively studied by researchers \citep{for, leb, lu}, and formulas for optimal shifts in \eqref{znumber} are given explicitly as
\begin{align}
  \label{optimalpq}
  p_j = M\left(\alpha \text{dn}\left(\frac{2j-1}{2k}K(\beta), \beta \right) \right), \quad
  q_j = M\left(-\alpha \text{dn}\left(\frac{2j-1}{2k}K(\beta), \beta \right) \right),
\end{align}
for $1 \leq j \leq k$, where $\alpha = -1 + 2\gamma + 2\sqrt{\gamma^2-\gamma}$ with $\gamma = \frac{|c-a||d-b|}{|c-b||d-a|}$, $\beta = \sqrt{1-1/\alpha^2}$, $K(\cdot)$ is the complete elliptic integral of the first kind, dn$(\cdot, \cdot)$ is the Jacobi elliptic function of the third kind, and $M$ is the M\"obius transformation that maps $\left\{-\alpha, -1, 1, \alpha\right\}$ to $\left\{a, b, c, d\right\}$. With optimal shifts in \eqref{optimalpq}, the upper bound for the Zolotarev number can be derived \citep{bec}
\begin{align}
  \label{znumberbound}
  Z_k([a, b], [c, d]) \leq 4\left[\exp \left(\frac{\pi^2}{2 \log \left(16 \gamma\right)}\right)\right]^{-2k}.
\end{align}
Therefore, if we expect a reduction of $\epsilon$ in the relative error \eqref{relerror}, the least number of iterations needed can be deduced from \eqref{znumberbound}, that is
\begin{align}
  \label{optimalk}
  k = \left\lceil \frac{\log \left(16 \gamma\right) \log (4/\epsilon)}{\pi^2} \right\rceil.
\end{align}

When solving the Poisson equation using the US method \eqref{poissonsylbc}, operators must be truncated to finite dimensions. For notational convenience, the same truncations are used for both $x$- and $y$-directions, though our method and analysis can easily adapt to different degrees of freedom in each direction. Let the $n \times \infty$ projection operator be $\mP_n = (I_n , \boldsymbol{0})$, and the generalized Sylvester equation we solve is
\begin{align}
  \label{truncations}
  A_1 \tilde{X} B_2 + A_2 \tilde{X} B_1 &= F \\
  \label{reducedtruncations}
  \text{or~~} A_2^{-1}A_1 \tilde{X} + \tilde{X} B_1 B_2^{-1} &= A_2^{-1} F B_2^{-1},
\end{align}
where
\begin{align}
  \label{poissoncoeffs}
  \begin{aligned}
    &A_1 = B_1^T = \mP_n \mS_1 \mS_0 \mT_{D} \mP_n^{\top}, \quad B_2^T = A_2 = \mP_n \mD_2 \mT_{D} \mP_n^{\top}, \\
    &F = \mP_n \mS_1 \mS_0 \mF (\mP_n \mS_1 \mS_0)^{\top}, \quad \tilde{X} = \mP_n \tilde{\mX} \mP_n^{\top}.
  \end{aligned}
\end{align}
For \eqref{reducedtruncations}, the ADI iteration \eqref{adiiter} then becomes
\begin{align}
  \label{uscoeffs}
  \begin{aligned}
    (A_2^{-1} A_1 - p_j I_n) \tilde{X}_{j-1/2} &= A_2^{-1} F B_2^{-1} - \tilde{X}_{j-1} (B_1 B_2^{-1} + p_j I_n),\\
    \tilde{X}_{j}(B_1 B_2^{-1} + q_j I_n) &= A_2^{-1} F B_2^{-1} - (A_2^{-1} A_1 - q_j I_n) \tilde{X}_{j-1/2}.
  \end{aligned}
\end{align}
To avoid explicitly forming $A_2^{-1} A_1$ and $B_1 B_2^{-1}$, we multiply both equations by $A_2$ and $B_2$ on the left and the right, respectively. Denoting $\hat{X}_{j-1/2} = \tilde{X}_{j-1/2} B_2$ and $\hat{X}_{j} = A_2 \tilde{X}_{j}$, we derive the ADI iteration for the generalized Sylvester equation \eqref{truncations}:
\begin{align}
  \label{gadiiter}
  \begin{aligned}
    (A_1 - p_j A_2) \hat{X}_{j-1/2} &= F - \hat{X}_{j-1} (B_1 + p_j B_2),\\
    \hat{X}_{j}(B_1 + q_j B_2) &= F - (A_1 - q_j A_2) \hat{X}_{j-1/2}.
  \end{aligned}
\end{align}

The prototype of our algorithm is summarized in \Cref{alg:prototype}. The generalized ADI for the generalized Sylvester equation is also proposed in \citet{kno}. Note that solutions of banded systems are needed in lines 5--6, which cost $\mathcal{O}(n^2)$ operations per iteration. An alternative simplification for \eqref{truncations} is $A = A_1^{-1} A_2$, $B = B_1 B_2^{-1}$, and $F = A_1^{-1} F B_1^{-1}$, which is mathematically equivalent to \eqref{uscoeffs}. However, this approach is deprecated due to numerical instability and reasons that will be explained later in \cref{subsec:bcs}.

\begin{algorithm}[!t]
  \caption{US-ADI method for solving the Poisson equation \eqref{poissonD}}
  \label{alg:prototype}
  \begin{unlist}
    \item[Inputs: ] $A_1, B_1, A_2, B_2, F \in \mathbb{R}^{n \times n}$ as given in \eqref{poissoncoeffs}, estimations for spectra $[a, b]$ and $[c, d]$ where $\sigma(A_2^{-1}A_1) \subset [a, b]$ and $\sigma(-B_1B_2^{-1}) \subset [c, d]$, and a tolerance $\epsilon$
    \item[Outputs: ] An approximate solution $X \in \mathbb{R}^{n \times n}$ for the coefficients in \eqref{usapprox}
  \end{unlist}
  \hrule
  \begin{algorithmic}[1]
  \State Compute the number of iterations $k$ in \eqref{optimalk} with $a, b, c, d$, and $\epsilon$.
  \State Compute optimal shifts $p_j$ and $q_j$ ($j=1,\dots,k$) using \eqref{optimalpq} with $a, b, c, d$, and $k$.
  \State Set initial iteration $\tilde{X}_0 = 0$ and $\hat{X}_0 = A_1 \tilde{X}_0$.
  \For {$j=1,\dots,k$}
    \State Solve $(A_1 - p_j A_2) \hat{X}_{j-1/2} = F - \hat{X}_{j-1} (B_1 + p_j B_2)$ for $\hat{X}_{j-1/2}$.
    \State Solve $\hat{X}_{j}(B_1 + q_j B_2) = F - (A_1 - q_j A_2) \hat{X}_{j-1/2}$ for $\hat{X}_{j}$.
  \EndFor
  \State Solve $A_2 \tilde{X} = \hat{X}_k$ for $\tilde{X}$.
  \State $X = \mP_{n}\mT_{D}\mP_{n}^{\top} \tilde{X} \mP_{n}\mT_{D}^{\top}\mP_{n}^{\top}$
  \end{algorithmic}
\end{algorithm}

The comparisons between our method and that in \citet{for} are worth mentioning. Firstly, the coefficient matrices in \citet{for} are symmetric and thus normal, so their relative error \eqref{relerror} is independent of $\kappa(U_A)$ and $\kappa(U_B)$. However, this satisfying result is brought about by a particular ultraspherical series, which cannot be extended to other problems. Although our coefficient matrices \eqref{uscoeffs} are nonnormal, the eigenvectors $U_A$ and $U_B$ are well-conditioned since they are good approximations to the orthogonal eigenvectors of the continuous problem \eqref{secondeigfunc} (see below). Numerical results show that $\kappa(U_A)$ and $\kappa(U_B)$ are at a reasonable level and can be ignored safely when determining iteration numbers from \eqref{optimalk}.

Secondly, our method is capable of solving many different problems besides \eqref{poissonD}. For the Poisson equation with Neumann boundary conditions, only $\mT_D$ in \eqref{poissonsylbc} and \eqref{poissoncoeffs} needs to be changed to a proper transformation operator $\mT_N$ (see \cref{subsec:bcs}). In contrast, the method in \citet{for} cannot solve this problem.

Thirdly, the bounds $a$, $b$, $c$, and $d$ in \eqref{optimalk} are critical for determining the number of ADI iterations. The smaller bounds we have, the fewer iterations are required for the same tolerance. Despite using different series, both methods construct approximations to the second-order differential operator with homogeneous Dirichlet conditions. We provide a detailed analysis and derive sharp bounds on eigenvalues in \cref{subsec:est}, which are smaller than those in \citet{for}, while only crude estimations from an algebraic point of view are shown in \citet[Appendix B]{for}. This is another advantage of our method since fewer iterations are required when solving \eqref{truncations}.

Lastly, their algorithm relies on Legendre polynomials, with transformations between Chebyshev and Legendre series each requiring $\mathcal{O}(n^2 (\log n)^2)$ operations. This cost overshadows the $\mathcal{O}(n^2 \log n)$ complexity of their solving process. In contrast, our method's transformation from the recombined basis to the canonical Chebyshev basis incurs only $\mathcal{O}(n^2)$ operations, preserving the overall $\mathcal{O}(n^2)$ complexity of our approach.

\subsection{Estimations for spectra}\label{subsec:est}
So far, it is assumed that the spectrum of the second-order US differential matrix is real and can be enclosed by an interval. Here, we provide an analysis and computable bounds for it. It is well-known that the eigenvalues of the continuous second-order derivative operator with homogeneous Dirichlet boundary conditions at $x=\pm 1$,
\begin{align}
  \label{secondeigproblem}
  \mL u(x) := u''(x) = \tilde{\lambda} u(x), \quad x \in [-1, 1], \quad u(\pm1) = 0,
\end{align}
are given by
\begin{align}
  \label{secondeig}
  \tilde{\lambda}_k = - \frac{k^2 \pi^2}{4}, \quad k = 1, 2, \dots,
\end{align}
with the corresponding normalized eigenfunctions
\begin{align}
  \label{secondeigfunc}
  u_k(x) = \begin{cases}
    \cos(k \pi x/2), & \text{for } k \text{ odd}, \\
    \sin(k \pi x/2), & \text{for } k \text{ even}.
  \end{cases}
\end{align}
While the eigenvalues of FD discretization of \eqref{secondeigproblem} on an $n$-grid are already known \citep{lev} and good estimations for those related to spectral collocation methods are given in \citet{wei}, only a loose lower bound for the spectrum of \eqref{secondeigproblem} approximated by the US method has been shown in \citet{che} before. Here, we derive sharp bounds for both ends following the ideas introduced in \citet{cha, wei}.

Since the second differential operator $\mD_2$ maps Chebyshev coefficients to ultraspherical $C^{(2)}$ coefficients, the eigenvalue problem \eqref{secondeigproblem} approximated by US method should be posed as
\begin{align}
  \label{ussecondeig}
  \begin{aligned}
      \left(\mL u_n(x), C^{(2)}_j(x)\right)_{C^{(2)}} &= \left(\lambda u_n(x), C^{(2)}_j(x)\right)_{C^{(2)}}, \quad j = 0, 1, \dots, n-2, \\
      u_n(\pm1) &= 0,
  \end{aligned}
\end{align}
where $u_n(x) \in \text{span}\left\{T_0(x), T_1(x), \dots, T_n(x)\right\}$ is an eigenfunction corresponding to eigenvalue $\lambda$, and $(\cdot, \cdot)_{C^{(2)}}$ is the weighted inner product with weight function $(1-x^2)^{3/2}$. Using the completeness of $C^{(2)}(x)$, we know that
\begin{align}
  \label{eigfunceq}
  \mL u_n(x) - \lambda u_n(x) = \xi C^{(2)}_{n-1}(x) + \eta C^{(2)}_{n}(x), \quad x \in [-1, 1],
\end{align}
since $\mL u_n(x) - \lambda u_n(x)$ is a polynomial of degree no more than $n$. Note that $\lambda$ is nonzero; otherwise, we can deduce that $u_n(x) \equiv 0$ from $\mL u_n(x) = 0$ and $u_n(\pm1) = 0$. The eigenfunction $u_n(x)$ for \eqref{ussecondeig} is formally given by
\begin{align}
  \label{eigfunc}
  \begin{aligned}
    u_n(x) &= (\mL - \lambda)^{-1}\left(\xi C^{(2)}_{n-1}(x) + \eta C^{(2)}_{n}(x)\right) \\
           &= -\xi \lambda^{-1} \mG\left(C^{(2)}_{n-1}(x)\right) - \eta \lambda^{-1} \mG\left(C^{(2)}_{n}(x)\right),
  \end{aligned}
\end{align}
where $\mG$ is the linear operator defined by
\begin{align}
  \label{Linv}
  \mG = \sum_{k=0}^{\left\lfloor n/2 \right\rfloor} \lambda^{-k} \mL^k.
\end{align}
The sum in $\mG$ is finite since $\mL^k\left(C^{(2)}_{n}(x)\right) \equiv 0$ for $k > n/2$. The zero Dirichlet conditions satisfied by $u_n$ yield
\begin{align}
  \label{pqeq}
  \xi \mG\left(C^{(2)}_{n-1}\right) + \eta \mG\left(C^{(2)}_{n}\right) = 0 \quad \text{at} \quad x = \pm 1.
\end{align}
The system in \eqref{pqeq} has nontrivial solutions for $\xi$ and $\eta$ if and only if
\begin{align}
  \label{pqdet}
  \mG\left(C^{(2)}_{n-1}\right)\Big|_{x=1}\mG\left(C^{(2)}_{n}\right)\Big|_{x=-1} - \mG\left(C^{(2)}_{n-1}\right)\Big|_{x=-1}\mG\left(C^{(2)}_{n}\right)\Big|_{x=1} = 0.
\end{align}
Assuming $n$ is even (odd $n$ can be analyzed similarly), $C^{(2)}_{n-1}$ is odd and $C^{(2)}_{n}$ is even, which implies that $\mG(C^{(2)}_{n-1})$ is odd and $\mG(C^{(2)}_{n})$ is even. Therefore, \eqref{pqdet} leads to $\mG\left(C^{(2)}_{n-1}\right)\Big|_{x=1}\mG\left(C^{(2)}_{n}\right)\Big|_{x=1} = 0$, that is,
\begin{align}
    \label{character}
    \sum_{k=0}^{n/2-1} a_{n/2-1-k} \lambda^k = 0 \quad \text{or} \quad \sum_{k=0}^{n/2} b_{n/2-k} \lambda^k = 0,
\end{align}
where
\begin{align}
    \label{abcoeff}
    \begin{aligned}
      a_k &= \mL^{k}\left(C^{(2)}_{n-1}\right)\Big|_{x=1} = \frac{2^{2k} \Gamma(n+2k+3) (2k+1)!}{\Gamma(4k+4) (n-2k-1)!}, \\
      b_k &= \mL^{k}\left(C^{(2)}_{n}\right)\Big|_{x=1} = \frac{2^{2k} \Gamma(n+2k+4) (2k+1)!}{\Gamma(4k+4) (n-2k)!},
    \end{aligned}
\end{align}
and $\Gamma(\cdot)$ is the Gamma function. All $n-1$ eigenvalues for \eqref{ussecondeig} are provided by the characteristic equations \eqref{character}, which can be proved to be real, negative, and distinct folllowing \citet{cha}. However, bounds on them are not provided in \citet{cha}. With the guarantee that eigenvalues are real and negative, we can employ Newton's root bound theorem on \eqref{character}, which states that the roots of the polynomial $\lambda^m + c_1 \lambda^{m-1} + \dots + c_m = 0$ satisfy
\begin{align}
  \label{Newtonbound}
  |\lambda| \leq \sqrt{\lambda_1^2 + \dots + \lambda_m^2} = \sqrt{c_1^2 - 2c_2},
\end{align}
provided that all of them are real. The upper bounds for all eigenvalues of \eqref{ussecondeig} can be obtained by applying \eqref{Newtonbound} to \eqref{character}, and the lower bounds by the same method after setting $\mu = \lambda^{-1}$. Since all eigenvalues are negative, we have
\begin{align}
  \label{secondNewton1}
  -\left(\frac{a_1^2}{a_0^2} - \frac{2a_2}{a_0}\right)^{1/2} \leq \lambda \leq -\left(\frac{a_{n/2-2}^2}{a_{n/2-1}^2} - \frac{2a_{n/2-3}}{a_{n/2-1}}\right)^{-1/2}, \\
  \label{secondNewton2}
  -\left(\frac{b_1^2}{b_0^2} - \frac{2b_2}{b_0}\right)^{1/2} \leq \lambda \leq -\left(\frac{b_{n/2-1}^2}{b_{n/2}^2} - \frac{2b_{n/2-2}}{b_{n/2}}\right)^{-1/2}.
\end{align}
Substituting \eqref{abcoeff} into \eqref{secondNewton1} and \eqref{secondNewton2} yields
\begin{align}
  \label{useigbound1}
  &\begin{aligned}
    &-\left(\frac{(n - 1)(n - 2)(n + 4)(n + 3)(29n^4 + 116n^3 + 1757n^2 + 3282n - 22824)}{121275}\right)^{1/2} \\
    &\qquad \leq \lambda \leq -\left(\frac{720(n^3 - n^2)}{8n^3 - 8n^2 - 45}\right)^{1/2},
  \end{aligned} \\
  \label{useigbound2}
  &\begin{aligned}
    &-\left(\frac{n(n - 1)(n + 5)(n + 4)(29n^4 + 232n^3 + 2279n^2 + 7260n - 17640)}{121275}\right)^{1/2} \\
    &\qquad \leq \lambda \leq -\left(\frac{48(n^3 + 2n^2 + n)}{8n^3 + 16n^2 + 8n - 3}\right)^{1/2}.
  \end{aligned}
\end{align}

\begin{figure}[!t]
  \centering
  \subfloat{\label{subfig:est}
  \includegraphics[width = 6cm]{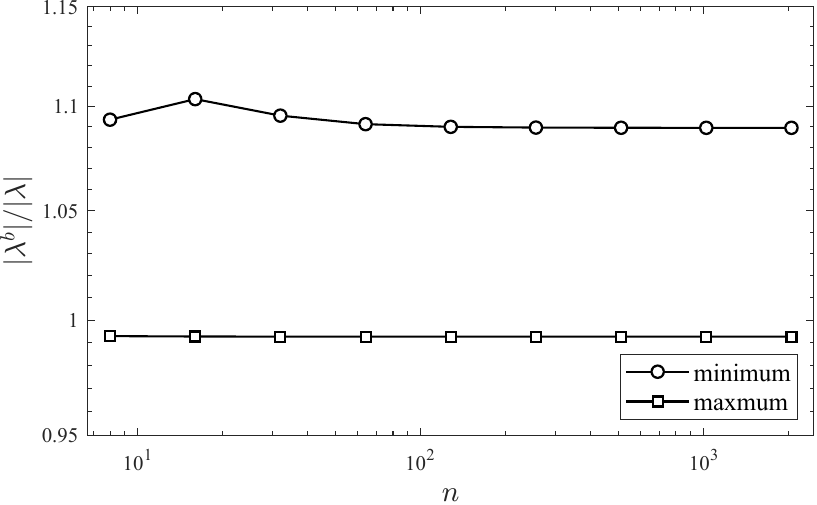}} \quad
  \subfloat{\label{subfig:convergence}
  \includegraphics[width = 6cm]{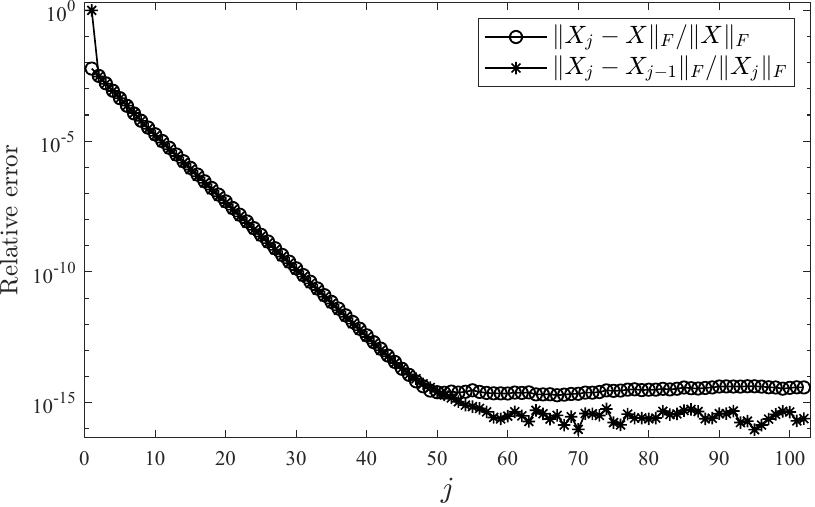}}
  \caption{Left: the ratios of bounds in \eqref{useigbound2} and extreme eigenvalues computed numerically for different finite truncations $n$. Right: relative true error and relative increment error during the ADI iterations when solving a Poisson equation with \Cref{alg:prototype} and sufficiently large $n$.}
  \label{fig:est}
\end{figure}

An immediate observation is that the range in \eqref{useigbound1} is contained within the range in \eqref{useigbound2}. Therefore, we take both ends in \eqref{useigbound2} as our bounds for the eigenvalues in \eqref{ussecondeig}, denoted by $\lambda^b_{\min}$ and $\lambda^b_{\max}$, respectively. \Cref{subfig:est} displays the ratios of the bounds in \eqref{useigbound2} and the extreme eigenvalues computed via numerical methods. The differences between the estimates and the true eigenvalues are less than 9 percent, implying that the bounds in \eqref{useigbound2} are sharp.

\begin{remark}
  In \citet{for}, bounds for the eigenvalues of their approximation to \eqref{secondeigproblem} are derived by a different method and are relatively loose. Consequently, empirical bounds are used in the practical implementation, specifically in \texttt{chebop2.poisson} in Chebfun \citep{dri}. Surprisingly, the analysis presented here can also be applied to their method, replacing $C^{(2)}$ in \eqref{eigfunceq} with $\tilde{C}^{(3/2)}$, the carefully chosen basis. Numerical computations show that the bounds derived by the analysis above match the empirical ones quite well. However, the spectrum of their method is still larger than ours, indicating that a larger value of $k$ in \eqref{optimalk} is required to achieve the same tolerance $\epsilon$.
\end{remark}

\begin{remark}
  \label{re:reversebound}
  Note that \cref{alg:prototype} is equivalent to solving a Sylvester equation with $A_2^{-1} A_1$ and $B_1 B_2^{-1}$ in \eqref{poissoncoeffs}, where these terms represent the inverses of the second-order differential operators. Consequently, we should use the reciprocals of the bounds in \eqref{useigbound2} to generate shifts in the ADI iterations. 
\end{remark}

\subsection{Optimal complexity}\label{subsec:opt}

From \eqref{useigbound2} and \eqref{optimalk}, the ADI method requires $k = \mathcal{O}(\log n \log (1/\epsilon))$ iterations to achieve a relative error \eqref{relerror} below a tolerance $\epsilon$, yielding an asymptotic complexity of $\mathcal{O}(n^2 \log n \log (1/\epsilon))$ for solving \eqref{truncations}\citep{for}. However, this estimate does not fully capture the behavior when solutions to PDEs exhibit sufficient smoothness. In \Cref{subfig:convergence}, we plot the relative true error and relative increment error for the Poisson equation $u_{xx} + u_{yy} = f$ solved via \Cref{alg:prototype} with Dirichlet boundary conditions. Here, $f$ is chosen such that the exact solution is $u(x, y) = 10 \exp(2x) \cos(2y)$, the target tolerance is double-precision $\epsilon \approx 2.2 \times 10^{-16}$, and the truncation dimension in \eqref{truncations} is $n = 2048$---far exceeding the resolution needed for the underlying solution. Using \Cref{alg:prototype}, \eqref{optimalk} predicts $k = 102$ iterations, yet the error curves stabilize after just 55 iterations. When repeated with the same tolerance but varying $n > 100$, the error consistently decreases monotonically over the first 55 iterations, plateauing at approximately $10^{-14}$ thereafter. This suggests that, for large $n$, a fixed number of iterations suffices for convergence under a smoothness assumption on $u$. We substantiate this observation below. For a related analysis of the Poisson equation using finite-difference (FD) discretization, see \citet{lyn}.

Suppose that real matrices $A$ and $-B^T$ have distinct eigenvalues $\lambda_i$ and $\mu_j$ and corresponding right eigenvectors $g_i$ and $h_j$, $i, j = 1, \dots, n$. It is easily verified that $g_i h_j^T$ form a complete basis for $\mathbb{R}^{n \times n}$. The initial error $E^0 := X - X_0$ can be expanded in terms of these eigenvectors as
\begin{align}
  \label{initialerror}
  E^0 = X - X_0 = \sum_{j=1}^{n} \sum_{i=1}^{n} e_{ij}^0 \left(g_i h_j^T\right).
\end{align}
By \eqref{adierror}, we have
\begin{align}
  \label{errork}
  \begin{aligned}
    E^k := X - X_k &= s_k(A) \left(\sum_{j=1}^{n} \sum_{i=1}^{n} e_{ij}^0 \left(g_i h_j^T\right)\right) s_k(-B)^{-1} \\
                   &= \sum_{j=1}^{n} \sum_{i=1}^{n} e_{ij}^0 s_k(\lambda_i) \left(g_i h_j^T\right) s_k(\mu_j)^{-1},
  \end{aligned}
\end{align}
where $s_k(\cdot)$ is defined in \eqref{adisk}. Without loss of generality, we assume that $g_i$ and $h_j$ are normalized, which leads to
\begin{align}
  \label{errorkbound}
  \left\| E^k \right\|^2 \leq \sum_{j=1}^{n} \sum_{i=1}^{n} \left(e_{ij}^0 s_k(\lambda_i)s_k(\mu_j)^{-1}\right)^2 \left\|g_i h_j^T\right\|^2 = \sum_{j=1}^{n} \sum_{i=1}^{n} \left(e_{ij}^0 \frac{s_k(\lambda_i)}{s_k(\mu_j)}\right)^2.
\end{align}
Note that \eqref{errorkbound} provides a more detailed bound than \eqref{relerror}, as it weights each term individually. Intuitively, for components with large initial errors $e_{ij}^0$, the reduction factor $s_k(\lambda_i) s_k(\mu_j)^{-1}$ should be smaller than for those with negligible errors. However, the analyses in \Cref{sec:adi} and \citet{for} overlook the influence of $e_{ij}^0$, treating all terms uniformly. Specifically, they aim to reduce every $s_k(\lambda_i) s_k(\mu_j)^{-1}$ below the tolerance, potentially leading to unnecessary iterations. For sufficiently smooth initial errors, the role of $e_{ij}^0$ becomes significant: error components with large $i$ or $j$ naturally fall below the tolerance, leaving only a subset dominant in the total error. This insight directs us to prioritize the substantial parts of the initial error. To quantify this smoothness, we introduce the following assumption.
\begin{assumption}
  \label{ass:smoothness}
  There are positive and bounded functions $\omega_1(i)$ and $\omega_2(j)$, both defined on $\{1, 2, \dots\}$, such that for any $n$
  \begin{align}
    \label{esmoothness}
    \left|e_{ij}^0\right| \leq \omega_1(i) \omega_2(j), \quad i=1,\dots,n, \quad j = 1, \dots, n.
  \end{align}
  Moreover, given $\epsilon$, there are positive integers $K_1(\epsilon)$ and $K_2(\epsilon)$ such that
  \begin{align}
    \label{wsmoothness}
    \sum_{i=K_1(\epsilon)}^{\infty} \omega_1(i)^2 < \epsilon, \quad \sum_{j=K_2(\epsilon)}^{\infty} \omega_2(j)^2 < \epsilon.
  \end{align}
\end{assumption}

An immediate corollary from \eqref{wsmoothness} is that there are constants $W_1$ and $W_2$, both independent of $\epsilon$, such that
\begin{align}
  \label{wbound}
  \sum_{i=1}^{\infty} \omega_1(i)^2 \leq W_1, \quad \sum_{j=1}^{\infty} \omega_2(j)^2 \leq W_2.
\end{align}
Note that, for fixed $i$ and $j$, the initial error $e_{ij}^0$ varies with $n$ but converges as $n \to \infty$. This follows from the eigenvectors $g_i$ and $h_j$ serving as effective approximations to \eqref{secondeigfunc} for sufficiently large $n$ \citep[Thm.2]{wei}, causing $e_{ij}^0$ to approach the Fourier coefficients of the initial error. If the initial error exhibits sufficient smoothness, e.g., possessing a second derivative of bounded variation, the Fourier series converges \citep{tre}, lending credence to \Cref{ass:smoothness}. Additionally, \citet{cha} establishes that the eigenvalues in \eqref{ussecondeig} from successive approximations interlace. Consequently, each eigenvalue increases with $n$, and the interval $[\lambda_{\min}^b(n_0), \lambda_{\max}^b(n_0)]$ for a fixed $n_0$ encloses the first $n_0$ eigenvalues of any approximation with $n > n_0$, i.e.,
\begin{align}
  \label{eigcond}
  \lambda_{\min}^b(n_0) \leq \lambda_{i} \leq \lambda_{\max}^b(n_0), \quad i = 1, \dots, n_0
\end{align}
In addition, we have $[a,b] = [-d,-c] \subset [-\infty, 0]$ and $q_j = -p_j < 0$ in \eqref{znumber} for the Poisson equation and thus
\begin{align}
  \label{shiftcond}
  |s_k(x)| \leq 1, \quad x<0, \text{ and } |s_k(x)^{-1}| \leq 1, \quad x > 0.
\end{align}
With these assumptions and observations, the following result can be asserted, which shows that ADI iterations with a fixed number of iterations are sufficient for achieving a given tolerance regardless of $n$.
\begin{theorem}
  \label{thm:smoothness}
  With \Cref{ass:smoothness} and \eqref{wbound}, given $0< \epsilon < W_1 + W_2$, there is a positive number $k_{\epsilon}$ (independent of $n$) such that $\left\|E^{k_\epsilon}\right\| \leq \sqrt{3}\epsilon$.
\end{theorem}

\begin{proof}
  For $\epsilon_W = \epsilon^2 / (W_1 + W_2)$, there are positive integers $\tilde{K}_1 = K_1(\epsilon_W)$ and $\tilde{K}_2 = K_2(\epsilon_W)$ such that the inequalities in \eqref{wsmoothness} hold. Denote $n_{K} = \max(\tilde{K}_1, \tilde{K}_2)$, $a = -d = \lambda_{\min}^b(n_K)$, and $b = -c = \lambda_{\max}^b(n_K)$. With a prescribed tolerance $\epsilon_{\text{ADI}} = \epsilon / \sqrt{W_1 W_2}$, we can compute shifts $p_j$ and $q_j$, $j = 1,\dots, k_{\epsilon}$ as in \eqref{optimalpq}, where $k_{\epsilon}$ is deduced from \eqref{optimalk} as
  \begin{align}
    \label{nepsilon}
    k_{\epsilon} = \left\lceil \frac{\log \left(16 |c-a||d-b|/(|c-b||d-a|)\right) \log (4\sqrt{W_1 W_2}/\epsilon)}{\pi^2} \right\rceil.
  \end{align}
  Note that $k_{\epsilon}$ is independent of $n$. By \eqref{znumberbound}, we have
  \begin{align}
    \label{lambdabound}
    \left|\frac{s_{k_{\epsilon}}(\lambda_i)}{s_{k_{\epsilon}}(\mu_j)}\right| \leq Z_{k_{\epsilon}}([a,b], [c,d]) \leq \epsilon / \sqrt{W_1 W_2}, \quad \text{for } \lambda_i \in [a, b], \quad \mu_i \in [c, d].
  \end{align}
  From \eqref{errorkbound}, we can split the upper bound of the error into four parts, namely,
  \begin{align*}
    \left\| E^{k_\epsilon} \right\|^2 \leq \sum_{j=1}^{n} \sum_{i=1}^{n} \left(e_{ij}^0 \frac{s_{k_{\epsilon}}(\lambda_i)}{s_{k_{\epsilon}}(\mu_j)}\right)^2 =: E_1 + E_2 + E_3 + E_4,
  \end{align*}
  where
  \begin{align*}
    \begin{aligned}
      &E_1 = \sum_{j=1}^{\tilde{K}_2} \sum_{i=1}^{\tilde{K}_1}\left(e_{ij}^0 \frac{s_{k_{\epsilon}}(\lambda_i)}{s_{k_{\epsilon}}(\mu_j)}\right)^2, \quad E_2 = \sum_{j=1}^{\tilde{K}_2} \sum_{i=\tilde{K}_1+1}^{n}\left(e_{ij}^0 \frac{s_{k_{\epsilon}}(\lambda_i)}{s_{k_{\epsilon}}(\mu_j)}\right)^2, \\
      &E_3 = \sum_{j=\tilde{K}_2+1}^{n} \sum_{i=1}^{\tilde{K}_1}\left(e_{ij}^0 \frac{s_{k_{\epsilon}}(\lambda_i)}{s_{k_{\epsilon}}(\mu_j)}\right)^2, \quad E_4 = \sum_{j=\tilde{K}_2+1}^{n} \sum_{i=\tilde{K}_1+1}^{n}\left(e_{ij}^0 \frac{s_{k_{\epsilon}}(\lambda_i)}{s_{k_{\epsilon}}(\mu_j)}\right)^2.
    \end{aligned}
  \end{align*}
  Note that some terms may be empty for different $n$. For $E_1$, given that $\lambda_i \in [a, b]$, $i=1, \dots, \tilde{K}_1$ and $\mu_j \in [c, d]$, $j = 1, \dots, \tilde{K}_2$, we have
  \begin{align*}
    E_1 \leq \frac{\epsilon^2}{W_1 W_2} \sum_{j=1}^{\tilde{K}_2} \omega_2(j)^2 \sum_{i=1}^{\tilde{K}_1} \omega_1(i)^2 \leq \frac{\epsilon^2}{W_1 W_2} W_2 W_1 = \epsilon^2.
  \end{align*}
  Since $\lambda_i < 0$ and $\mu_j>0$ for all $i$ and $j$, we conclude from \eqref{shiftcond} that for $E_2$, $E_3$, and $E_4$, $\left(s_{k_{\epsilon}}(\lambda_i)s_{k_{\epsilon}}(\mu_j)^{-1}\right)^2 \leq 1$,
  which results in
  \begin{align*}
    E_2 &\leq \sum_{j=1}^{\tilde{K}_2} \omega_2(j)^2  \sum_{i=\tilde{K}_1+1}^{n} \omega_1(i)^2 \leq \frac{W_2\epsilon^2}{W_1 + W_2}, \quad E_3 \leq \frac{W_1\epsilon^2}{W_1 + W_2}, \\
    E_4 &\leq \sum_{j=\tilde{K}_2+1}^{n} \omega_2(j)^2 \sum_{i=\tilde{K}_1+1}^{n} \omega_1(i)^2 \leq \frac{\epsilon^4}{(W_1 + W_2)^2} \leq \epsilon^2,
  \end{align*}
  where \eqref{wsmoothness} corresponding to $\epsilon_W = \epsilon^2 / (W_1 + W_2)$ and $\epsilon < W_1 + W_2$ are used. Taken together, we have
  \begin{align}
    \label{enlarge}
    \| E^{k_{\epsilon}} \|^2 \leq \epsilon^2 + \frac{W_2\epsilon^2}{W_1 + W_2} + \frac{W_1\epsilon^2}{W_1 + W_2} + \epsilon^2 = 3\epsilon^2.
  \end{align}
  Therefore, the theorem holds no matter what $n$ is used for truncations.
\end{proof}

\Cref{thm:smoothness} elucidates the plateau observed in \Cref{subfig:convergence}, indicating that a fixed number of iterations suffices for the ADI method to solve PDEs with sufficiently smooth solutions. Consequently, the asymptotic complexity of \Cref{alg:prototype} reduces to $\mathcal{O}(n^2)$, which is optimal. However, \Cref{thm:smoothness} provides no practical method to determine $k_{\epsilon}$. Since the solution $u$, its coefficient matrix $X$, and the parameters $\omega_1$ and $\omega_2$ in \eqref{esmoothness} are unknown a priori, estimating $k_{\epsilon}$ directly is infeasible. To prevent oversolving, \Cref{sec:details} proposes a strategy for terminating ADI iterations efficiently. Notably, \Cref{thm:smoothness} also applies to the methods in \citet{for}.

\section{Implementation details}\label{sec:details}
Some useful points are introduced here to enhance the performance of our algorithm. The full algorithm for solving the Poisson equation automatically is summarized in \Cref{alg:auto}.

\begin{algorithm}[!t]
  \caption{Automatic algorithm for solving Poisson equations}
  \label{alg:auto}
  \begin{unlist}
    \item[Inputs: ] Operators $\mD_{i}$, $\mS_i$, $\mT$, and $\mP$, a tolerance $\epsilon$, and a positive integer $\tau$ for checking the relative increment error
    \item[Output: ] Approximate solution $u(x, y)$
  \end{unlist}
  \hrule
  \begin{algorithmic}[1]
    \State Initialize truncation dimension $n = 16$ and initial iterate $X_0 = 0$.
    \State Compute matrices $A_1$, $B_1$, $A_2$, $B_2$, and $F$ from \eqref{uscoeffs}.
    \State Solve for $X$ using \Cref{alg:prototype} with inputs $A_1$, $B_1$, $A_2$, $B_2$, $F$, $\epsilon$, and $X_0$, checking relative increment error every $\tau$ steps.
    \While{$u(x, y) = \sum_{i,j} X_{ij} T_i(y) T_j(x)$ remains unresolved}
      \State Update $n = 2n$ and recompute matrices in \eqref{truncations}.
      \State Generate shifts $p_j$ and $q_j$ as in \eqref{optimalpq}, then reverse their order.
      \State Set the initial iteration $X_0$ as the unresolved $X$ padded with zeros
      \State Solve for $X$ using \Cref{alg:prototype} with updated matrices, reversed shifts, $\epsilon$, and $X_0$, checking relative increment error every $\tau$ steps.
    \EndWhile
    \State \Return $u(x, y) = \sum_{i,j} X_{ij} T_i(y) T_j(x)$.
  \end{algorithmic}
\end{algorithm}

\subsection{Error check}\label{subsec:check}
The blind application of \Cref{alg:prototype} for solving PDEs would result in an $\mathcal{O}(n^2 \log n)$ complexity, whereas \Cref{thm:smoothness} ensures that a fixed number of iterations is sufficient, thus requiring only $\mathcal{O}(n^2)$ operations. To avoid unnecessary iterations, as shown in \Cref{subfig:convergence}, we propose monitoring the convergence of the relative increment error $\| X_{j} - X_{j-1} \| / \| X_{j} \|$ during the iterations since it is a reliable substitute for the relative true error $\| X - X_{j-1} \| / \| X \|$. Note that by \eqref{adierror}

\begin{align*}
  \begin{aligned}
    &X_j - X_{j-1} = X - X_{j-1} - (X - X_j)\\
    =&s_{j-1}(A) X s_{j-1}(-B)^{-1} - s_{j}(A) X s_{j}(-B)^{-1}\\
    =&s_{j-1}(A) (X - (A - q_j I_n)(A - p_j I_n)^{-1}X(B + p_j I_n)(B + q_j I_n)^{-1}) s_{j-1}(-B)^{-1}\\
    =&s_{j-1}(A) ((-p_j + q_j)(A - p_j I_n)^{-1}F(B + q_j I_n)^{-1}) s_{j-1}(-B)^{-1}\\
    =& \frac{-p_{j} + q_{j}}{-p_{j-1} + q_{j-1}} (A- q_{j-1} I_n)(A - p_{j} I_n)^{-1} (X_{j-1} - X_{j-2}) (B + p_{j-1} I_n)(B + q_{j} I_n)^{-1}\\
    =& c_{pq} (A- q_{j-1} I_n)(A - (p_{j-1} + \Delta p_{j}) I_n)^{-1} (X_{j-1} - X_{j-2}) (B + p_{j-1} I_n)(B + (q_{j-1} + \Delta q_{j}) I_n)^{-1},
  \end{aligned}
\end{align*}
where $\Delta p_{j} = p_{j} - p_{j-1}$, $\Delta q_{j} = q_{j} - q_{j-1}$ and $c_{pq} = (1 + \frac{-\Delta p_{j} + \Delta q_{j}}{-p_{j-1}+q_{j-1}})$. This implies that the increment error has almost the same rate of convergence as the true error provided that $\Delta p_{j}$ and $\Delta q_{j}$ are relatively small compared with $p_{j-1}$ and $q_{j-1}$, which is always true since the shifts in \eqref{optimalpq} are almost logarithmically equispaced.

Once $\| X_{j} - X_{j-1} \| / \| X_{j} \|$ is less than the prescribed tolerance or stagnates, the iterations are stopped. The additional cost is only $\mathcal{O}(n^2)$ per check. Since the number of iteration $k$ is usually within 100, we recommend checking the error every 5 or 10 steps for better performance of our algorithm, as suggested in \citet{li}. See lines 3 and 8 in \Cref{alg:auto}.

\subsection{Warm restart}\label{subsec:restart}
An optimal $n$ for the resolution of $u$ is not known in advance, and we have to solve \eqref{truncations} with progressively larger $n$. That is, we solve the problem with $n$-truncations of operators, check the solution $u$, and increase $n$ if necessary. An additional advantage of solving \eqref{truncations} using the ADI method is that a reasonable initial iteration for \eqref{adiiter} with larger $n$ can be easily constructed. Specifically, the last iteration of \Cref{alg:prototype} with smaller $n$, padded with zeros, can be used as the initial iteration for ADI method with larger $n$. See line 7 in \Cref{alg:auto}.

\subsection{Order of Shifts}\label{subsec:order}
From \eqref{errorkbound} and the proof of \Cref{thm:smoothness}, rapid convergence requires minimizing $s_k(\lambda_i) s_k(\mu_j)^{-1}$ first for eigenvalues $\lambda_i$ and $\mu_j$ associated with large initial error components $e_{ij}^0$. For most functions, low-frequency coefficients dominate. Thus, when starting with a zero initial iterate, we apply shifts in \eqref{gadiiter} corresponding to smaller-magnitude eigenvalues first, namely, $p_1$ and $q_1$ per \eqref{optimalpq} and \Cref{re:reversebound}. This ascending order is also endorsed by \citet{lyn}.

However, with a warm restart (see \Cref{subsec:restart}), the strategy changes. When transitioning to a larger truncation of \eqref{truncations} with a warm restart, the initial error predominantly comprises high-frequency terms, as low-frequency components are already diminished from the smaller Sylvester equation. To accelerate error reduction, we apply shifts in descending order, starting with $p_k$ and $q_k$ (see line 6 of \Cref{alg:auto}).

\subsection{Low-Rank Forms}\label{subsec:lowrank}
It is well established that the solution $X$ of \eqref{modelsyl} exhibits rapid singular value decay when the right-hand side $F$ and coefficient matrices have low-rank off-diagonal blocks \citep{mas}, particularly when the rank of $F$, denoted $r_F$, is significantly smaller than $n$ (i.e., $k r_F \ll n$). Indeed, methods exploiting low-rank forms of $X$ have proven highly efficient \citep{ben, li}. Our algorithm for solving \eqref{truncations} can be adapted to leverage a low-rank factorization $F = U_F V_F^T$, further reducing computational cost. Adopting the notation of factored ADI (fADI) from \citet{ben}, we express
\begin{align}
  \label{fadi}
  \begin{aligned}
    X_k &= Z_k D_k Y_k^T, ~D_k = \text{diag}((q_1 - p_1)I_{r_F}, (q_2 - p_2)I_{r_F}, \dots, (q_k - p_k)I_{r_F})\\
    Z_k &= (Z^{(1)} ~ Z^{(2)} ~ \cdots ~ Z^{(k)}), ~ Y_k = (Y^{(1)} ~ Y^{(2)} ~ \cdots ~ Y^{(k)}),
  \end{aligned}
\end{align}
where
\begin{align}
  \label{fadiiter1}
  \begin{aligned}
      Z^{(1)} &= (A - p_1 I_n)^{-1} U_F,~ Z^{(j+1)} = (A - p_{j+1} I_n)^{-1} (A - q_{j} I_n) Z^{(j)},\\
      Y^{(1)} &= (B + q_1 I_n)^{-T} V_F, ~ Y^{(j+1)} = (B + q_{j+1} I_n)^{-T} (B + p_{j} I_n)^T Y^{(j)},\\
  \end{aligned}
\end{align}
or
\begin{align}
  \label{fadiiter2}
  \begin{aligned}
      Z^{(1)} &= (A - p_1 I_n)^{-1} U_F,~ Z^{(j+1)} = Z^{(j)} + (p_{j+1} - q_j)(A - p_{j+1} I_n)^{-1} Z^{(j)},\\
      Y^{(1)} &= (B + q_1 I_n)^{-T} V_F, ~ Y^{(j+1)} = Y^{(j)} + (p_{j} - q_{j+1})(B + q_{j+1} I_n)^{-T} Y^{(j)}.\\
  \end{aligned}
\end{align}
Immediately, one can see that dominating operations needed for each iteration reduce from two matrix-matrix multiplications and two matrix-matrix solves in \cref{adiiter} to $2r_F$ matrix-vector solves in \cref{fadiiter2}. The error check can still be used since $X_{j} - X_{j-1} = (q_j - p_j)Z^{(j)} (Y^{(j)})^T$ can be easily computed. We note that fADI would not be favored if no low rank factorization of $F$ is available or $k r_F \geq n$.

\section{Extensions}\label{sec:more}

The US method introduced in \cref{sec:us} allows us to solve other problems besides the Poisson equation with homogeneous Dirichlet conditions \eqref{poissonD}. We show several important cases here to demonstrate the generality of our method.

\subsection{Poisson equation with Neumann and Robin conditions}\label{subsec:bcs}
To manage different boundary conditions, only transformation operator $\mT_D$ in \eqref{poissonsylbc} should be changed. If left Dirichlet and right Neumann conditions are imposed, the elements of the matrix in \eqref{reducedBC} are $\mB_{1i} = (-1)^i, ~ \mB_{2i} = i^2, ~ i = 0, 1, \dots,$ while for left Dirichlet and right Robin conditions, namely $u(-1) = u(1) + \theta u'(1) = 0$, the coefficients would be $\mB_{1i} = (-1)^i, ~ \mB_{2i} = 1 + \theta i^2, ~ i = 0, 1, \dots.$
To get a well-conditioned transformation operator, we fix $v_{k+2}^k = 1/(2k+4), \, k=0,1,\dots$ in \eqref{reducedBC} (which are the same as elements of the preconditioner proposed in \citet[\S 4]{olv}, implying that the effects of the preconditioner are combined in the transformation operator with this choice of $v_{k+2}$), and the transformation operator has the form
\begin{align}
  \label{transNR}
  \mT_{N} = \begin{pmatrix}
    t^N_0 \vphantom{\ddots} & & \\
    s^N_0 \vphantom{\ddots} & t^N_1 & \\
    r_0 & s^N_1 & \ddots \\
     & r_1 & \ddots \\
     &  & \ddots \\
  \end{pmatrix}, \quad
  \mT_{R} = \begin{pmatrix}
    t^R_0 \vphantom{\ddots} & & \\
    s^R_0 \vphantom{\ddots} & t^R_1 & \\
    r_0 & s^R_1 & \ddots \\
     & r_1 & \ddots \\
     &  & \ddots \\
  \end{pmatrix},
\end{align}
where 
\begin{align*}
  t^N_k &= -\frac{2k^2 + 6k + 5}{2(k + 2)(2k^2 + 2k + 1)}, \quad s^N_k = -\frac{2(k + 1)}{(k + 2)(2k^2 + 2k + 1)}, \quad r_k = \frac{1}{2k+4},\\
  t^R_k &= -\frac{2\theta k^2 + 6\theta k + 5\theta + 2}{2(k + 2)(2\theta k^2 + 2\theta k + \theta + 2)}, \quad s^R_k = -\frac{2(\theta + \theta k)}{(k + 2)(2\theta k^2 + 2\theta k + \theta + 2)},
\end{align*}
for $k=0,1,\dots$. Note that the elements on the main diagonal of both $\mD_2 \mT_N$ and $\mD_2 \mT_R$ are all 1s, ensuring the well-conditioness of our method.

The spectral estimates for coefficient matrices in \Cref{subsec:est} require adjustment to accommodate general boundary conditions. For example, with left Dirichlet and right Neumann conditions, the equations in \eqref{pqeq} become:
\begin{align*}
  \xi \mG\left(C^{(2)}_{N-1}\right)\Big|_{x=-1} + \eta \mG\left(C^{(2)}_{N}\right)\Big|_{x=-1} = \xi \mG'\left(C^{(2)}_{N-1}\right)\Big|_{x=1} + \eta \mG'\left(C^{(2)}_{N}\right)\Big|_{x=1} = 0
\end{align*}
The computation of characteristic polynomials and their coefficients remains valid. However, these general boundary conditions do not guarantee real eigenvalues. Numerical results indicate that some interior eigenvalues (excluding the extremes) appear as complex conjugates with significant imaginary parts. Nevertheless, the bounds in \eqref{Newtonbound} remain computable, as the extreme eigenvalues are real and dominate the square sum in \eqref{Newtonbound}. These complex eigenvalues challenge the use of optimal shifts from \eqref{optimalpq}, which assume real spectra for the Sylvester equation's coefficient matrices. Fortunately, we can mitigate the influence of complex numbers by taking the inverse of the operators to make their eigenvalues almost real. According to \citet{wac2}, we can safely use the optimal shifts \eqref{optimalpq} if the spectra are closely around the real axis. This is exactly what we do in \eqref{uscoeffs} as $A = A_2^{-1}A_1$ stands for the inverse of the one-dimensional Poisson operator.

\subsection{Poisson equation with separable coefficients}\label{subsec:variable}

As established in Sections~\ref{sec:adi} and~\ref{subsec:bcs}, the alternating direction implicit (ADI) method requires that the operators in the $x$- and $y$-directions have separate and nearly real spectra to ensure solvability. Our method extends this framework to a class of equations with non-constant, separable coefficients, such as:
\begin{align}
  \label{separablecoeffs}
  u_{xx} + u_{yy} - \bigl(\rho_1(x) + \rho_2(y)\bigr) u = f, \quad u(\pm 1, \cdot) = u(\cdot, \pm 1) = 0.
\end{align}
Since the variable coefficients of the zero-order term are separable, we can define two directional operators: $u_{xx} - \rho_1(x)u$ for the $x$-direction and $u_{yy} - \rho_2(y)u$ for the $y$-direction. For a wide range of functions $\rho_1$ and $\rho_2$, the spectra of these operators remain disjoint. For instance, if $\rho_1(x) = x^2$ and $\rho_2(y) = y^4$, the eigenvalues of $-\rho_1(x)u$ and $-\rho_2(y)u$ lie within $[-1, 0]$, ensuring that the spectra of both one-dimensional operators fall within the interval $[\lambda^b_{\min} - 1, \lambda^b_{\max}]$ \citep{gol2}. 

To accommodate these variable coefficients, the operator $\mD_2$ in \eqref{poissoncoeffs} must be modified for each direction: $\mD_2 - \mS_1 \mS_0 \mathcal{M}_0[\rho_1]$ for the $x$-direction and $\mD_2 - \mS_1 \mS_0 \mathcal{M}_0[\rho_2]$ for the $y$-direction, where $\mathcal{M}_0$ denotes the multiplication operator \citep{olv, qin}. Notably, the screened Poisson equation \citep{kno} is a special case of \eqref{separablecoeffs} with constant coefficients $\rho_1(x) = \rho_2(y) = -\omega^2 / 2$.

\subsection{Fourth-order equation}\label{subsec:fourth}

The fourth-order equation with Dirichlet and Neumann conditions is posed as
\begin{align}
  \label{fourth}
  \begin{aligned}
    u_{xxxx} + u_{yyyy} &= f, \quad (x, y) \in [-1, 1]^2,\\
    u(\pm 1, \cdot) = u(\cdot, \pm 1) &= 0, \quad \frac{\partial u}{\partial x}(\pm 1, \cdot) = \frac{\partial u}{\partial y}(\cdot, \pm 1) = 0.
  \end{aligned}
\end{align}
For a recombined basis satisfying the restrictions in \eqref{fourth}, we may take \citep{she2}
\begin{align}
  \label{fourthbasis}
  \psi^k(x) = T_k(x) - \frac{2(k+2)}{k+3}T_{k+2} + \frac{k+1}{k+3}T_{k+4}, \quad k = 0, 1, \dots,
\end{align}
which is also the choice of our procedure \eqref{reducedBC} with fixed $v^k_k=1$. With the transformation operator $\mT_F$ constructed from \eqref{fourthbasis}, differential operator $\mD_4$, and conversion operator $\mS_i$ for $i = 0, 1, 2, 3$, the fourth-order equation \eqref{fourth} is equivalent to
\begin{align}
  \label{fourthsylbc}
  \mS_3 \mS_2 \mS_1 \mS_0 \mT_{F} \tilde{\mX} (\mD_4 \mT_{F})^{\top} + \mD_4 \mT_{F} \tilde{\mX} (\mS_3 \mS_2 \mS_1 \mS_0 \mT_{F})^{\top} = \mS_3 \mS_2 \mS_1 \mS_0 \mF (\mS_3 \mS_2 \mS_1 \mS_0)^{\top},
\end{align}
where $\mX = \mT_F \tilde{\mX} \mT_{F}^{\top}$ is the coefficients of $u$ in tensor product Chebyshev basis. After the estimations for the spectrum of the inverse of the fourth-order differential operator with Dirichlet and Neumann conditions are obtained as in \cref{subsec:est}, we can solve for finite truncations of $\mX$ using \Cref{alg:prototype}.

\subsection{Nonhomogeneous boundary conditions}\label{subsec:nonhomogeneous}

For problems with nonhomogeneous boundary conditions, the system can be transformed into an equivalent problem with homogeneous boundary conditions by subtracting an interpolation of the boundary data from the true solution \citep{for, she1}. After obtaining the approximate solution to the modified equation with homogeneous boundary conditions and the adjusted right-hand side, the true solution is recovered by adding back the interpolation. This process incurs a computational cost of $\mathcal{O}(1)$ for constructing the interpolation and the associated multiplications and subtractions, which does not significantly impact the overall efficiency of our algorithm.

\section{Numerical results}\label{sec:exp}

We demonstrate the performance of our fast Poisson solver through four numerical examples. All experiments are conducted in Julia v1.10.2 on a desktop equipped with a 4.1 GHz Intel Core i7 CPU and 16 GB of RAM. The proposed approach, which combines the US method with basis recombination (\eqref{poissonsylbc} and \eqref{fourthsylbc}), is referred to as the ``new method'' or simply ``new''. Unless otherwise specified, the default solver for the new method is \Cref{alg:prototype} and LAPACK routines \texttt{gbtrf} and \texttt{gbtrs} are used for banded matrix factorization and solution. For comparison, we also implement the solvers from \citet{for, tow} in Julia, denoted as FT and TO, respectively. For error evaluation, we use a tolerance of $\tau = 10$ across all examples. Execution times are measured using the \texttt{BenchmarkTools.jl} package. The source code for the new method is publicly available at \url{https://github.com/ouyuanq/optimalPoisson}.

\subsection{Optimality}

\begin{figure}[!t]
  \centering
  \subfloat{\label{subfig:ex1time}
  \includegraphics[width = 6.5cm]{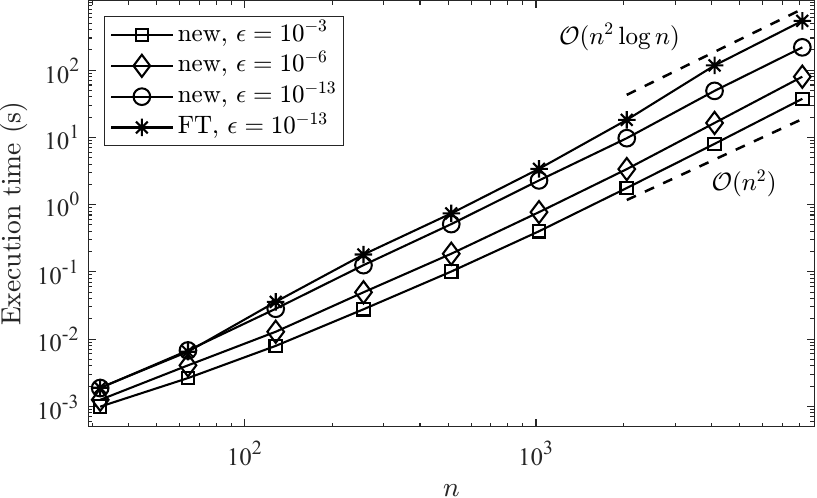}} \quad
  \subfloat{\label{subfig:ex1accuracy}
  \includegraphics[width = 6.5cm]{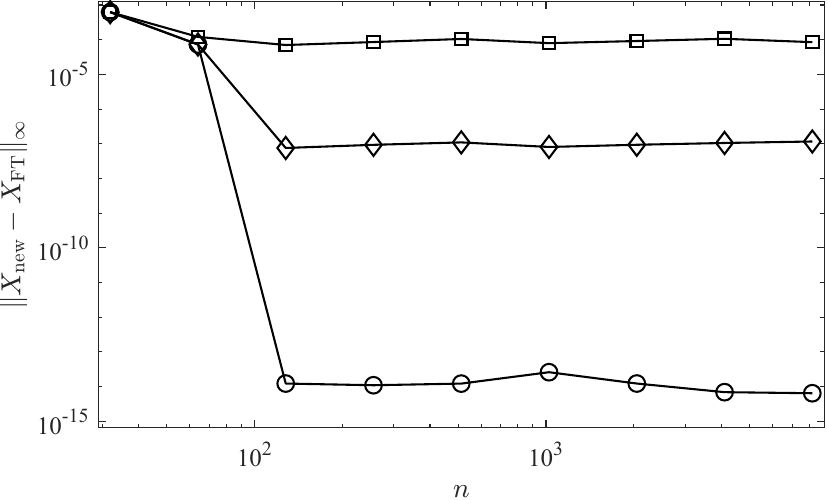}}
  \caption{Left: execution times for solving \eqref{ex1} using new method with different tolerances $\epsilon$ and the FT method with the tightest tolerance. Right: differences between solutions of new method and that of the FT method.}
  \label{fig:ex1}
\end{figure}

The first example is a Poisson equation with homogeneous Dirichlet boundary conditions, taken from \citet{for}:
\begin{align}
  \label{ex1}
  u_{xx} + u_{yy} = -100x \sin(20 \pi x^2 y) \cos(4 \pi (x + y)), \quad u(\pm 1, \cdot) = u(\cdot, \pm 1) = 0.
\end{align}
We solve this problem using the new method with varying truncation sizes and tolerances. For comparison, we also apply the FT method \citep{for} with its tightest tolerance settings. The execution times and the differences between the solutions obtained from the two methods are presented in \Cref{fig:ex1}. The left panel clearly illustrates the $\mathcal{O}(n^2)$ complexity of the new method. Notably, problems with millions of unknowns can be solved in a minute. In contrast, the FT method exhibits a complexity of $\mathcal{O}(n^2 \log n)$, consistent with the analysis in \citet{for}. This advantage in the new method arises from  the error-checking mechanism described in \Cref{subsec:check}, which eliminates unnecessary iterations. With the tightest tolerance $\epsilon=10^{-13}$, the FT method requires 34 to 105 iterations, while the new method needs 33 to 103. When $n$ is small, the differences in execution time between the new method and FT method are indistinguishable. However, the technique in \Cref{subsec:check} terminates ADI iterations at the 60th step for $n > 256$, reducing total iterations by up to one-third. While this technique could also enhance the FT solver, we exclude the cost of transforming to the Chebyshev basis here: $\mathcal{O}(n^2)$ for the new method (cf. line 9 of \Cref{alg:prototype}) versus $\mathcal{O}(n^2 (\log n)^2)$ for the FT method. Although this cost is minor relative to the solving phase for the values of $n$ considered, the transformation back to Chebyshev coefficients dominates the FT method's complexity breakdown and poses nontrivial implementation challenges. The new method circumvents these issues, requiring only two banded matrix multiplications. The right panel of \Cref{fig:ex1} shows that the differences between our solutions and those of the FT method remain below the prescribed tolerance $\epsilon$, confirming that the new method matches the FT method's accuracy.

\subsection{Adaptivity to various boundary conditions}

\begin{figure}[!t]
  \centering
  \subfloat{\label{subfig:ex2time}
  \includegraphics[width = 6.5cm]{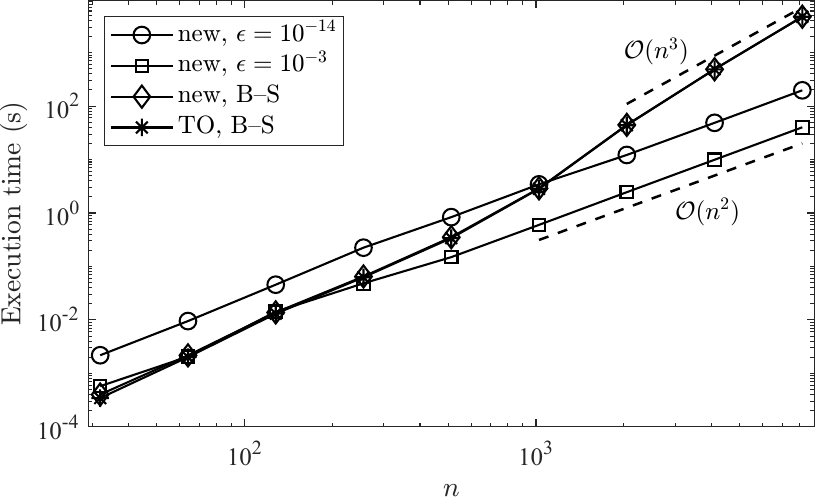}} \quad
  \subfloat{\label{subfig:ex2accuracy}
  \includegraphics[width = 6.5cm]{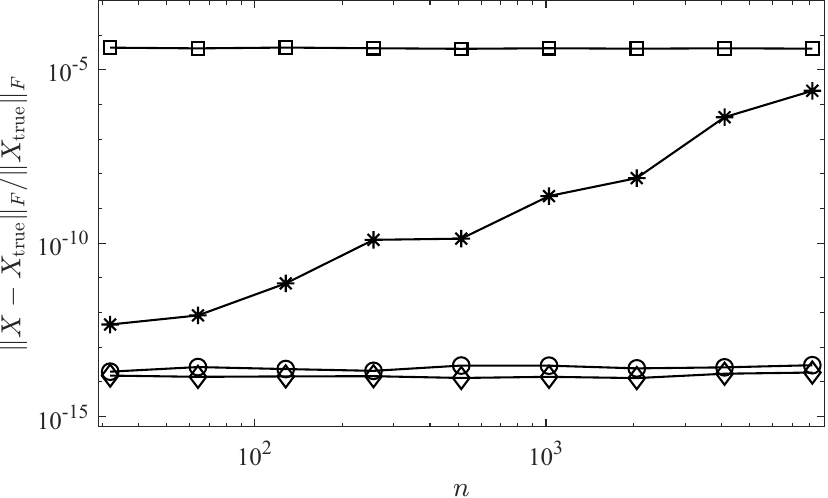}}
  \caption{Left: execution times for solving \eqref{ex2} using new method with different tolerances and the B--S algorithm with new method and TO method. Right: relative errors of different methods.}
  \label{fig:ex2}
\end{figure}

Next, we consider a Poisson equation with mixed boundary conditions:
\begin{align}
  \label{ex2}
  \begin{aligned}
    &u_{xx} + u_{yy} = 0, \quad u(-1, y) = 10 e^{-2} \cos(2y), \quad \frac{\partial u}{\partial x}(1, y) = 20 e^{2} \cos(2y), \\
    &u(x, -1) = 10 \cos(-2) e^{2x}, \quad u(x, 1) + \frac{\partial u}{\partial y}(x, 1) = (10 \cos(2) - 20 \sin(2)) e^{2x}.
  \end{aligned}
\end{align}
This problem features Dirichlet, Neumann, and Robin boundary conditions, designed to yield the exact solution $u(x, y) = 10 e^{2x} \cos(2y)$. Since the FT method cannot handle \eqref{ex2}, we solve it using the new method and the TO method. In addition to the default ADI solver, we employ the Bartels--Stewart (B--S) algorithm, invoked via Julia's \texttt{sylvester} function, to solve the truncated system in \eqref{reducedtruncations} for both the new method and TO method. Execution times and accuracy results are reported in \Cref{fig:ex2}. The left panel's last two curves, which overlap, reflect the B--S algorithm’s $\mathcal{O}(n^3)$ complexity, as it treats all coefficient matrices as full, regardless of their structure. In contrast, the new method achieves an optimal $\mathcal{O}(n^2)$ complexity. However, for small $n$ and tight tolerances, the new method is less efficient, though it can still provide a rough approximate solution with reduced computational effort. At a tolerance of $10^{-3}$, the ADI solver's solution time becomes comparable to that of the B--S algorithm, even for small to medium truncation sizes.

The right panel of \Cref{fig:ex2} shows that the ADI solver rivals the B--S algorithm in accuracy for our discretization. Clearly, the TO method exhibits a loss in relative error, which worsens as $n$ increases. For the new method, the accuracy is not altered no matter which solver is used. This indicates that the basis recombination, introduced in \Cref{sec:us} and \Cref{subsec:bcs}, is better conditioned than the TO method for Neumann and Robin boundary conditions. We posit that the approach in \Cref{subsec:bcs} can be extended to any linear boundary conditions while preserving the well-conditioned nature of the US method. In summary, the new method excels in both speed and accuracy.

\subsection{Weak singularity}

\begin{figure}[!t]
  \centering
  \subfloat{\label{subfig:ex3relerr}
  \includegraphics[width = 6.5cm]{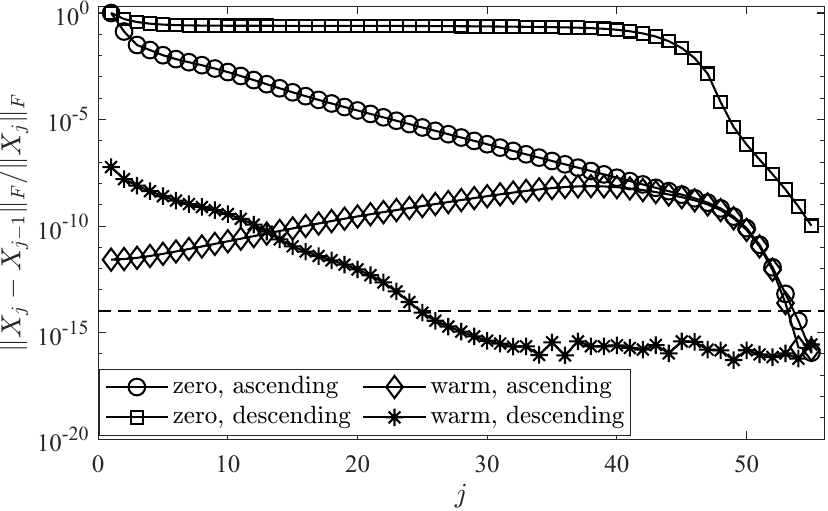}} \quad
  \subfloat{\label{subfig:ex3wr}
  \includegraphics[width = 6.5cm]{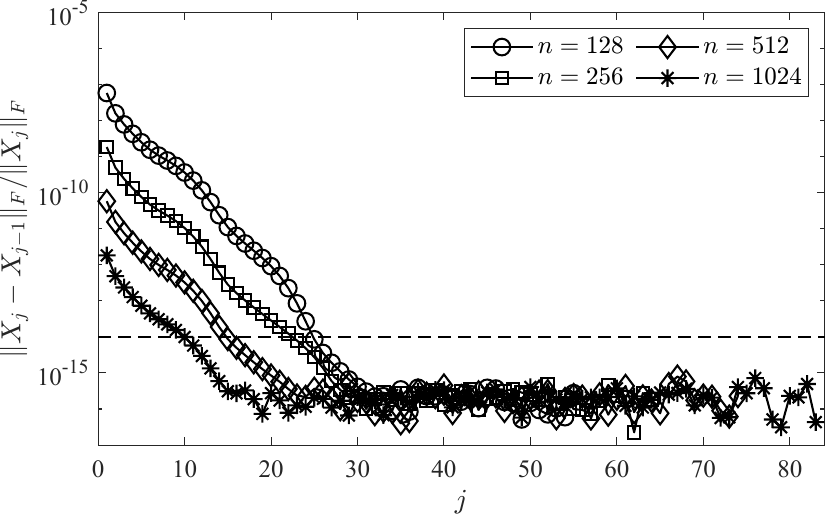}}
  \caption{Left: relative increment errors of ADI iterations for different initial iterations and different orders of application of shifts. Right: relative increment errors in the solution using warm restart and descending order of shifts.}
  \label{fig:ex3}
\end{figure}

For the third example, we consider a Poisson equation with separable coefficients and weak corner singularities:
\begin{align}
  \label{ex3}
  u_{xx} + u_{yy} + (-100 x^2 + \cos(\pi y)) u = 1, \quad u(\pm 1, \cdot) = u(\cdot, \pm 1) = 0.
\end{align}
We solve \eqref{ex3} using the techniques from \Cref{subsec:variable} and \Cref{alg:auto}, with a tolerance of $\epsilon = 10^{-14}$ and no error checking enforced. Relative errors are recorded with and without the methods from \Cref{subsec:restart} and \Cref{subsec:order}. The initial truncation size is $n = 16$, doubled iteratively until a resolved solution is achieved at $n = 1024$. For $n > 16$, we solve the subsequent Sylvester equation \eqref{truncations} in four ways: zero initial iteration with shifts in ascending or descending order, and warm restart with shifts in ascending or descending order. The first system ($n = 16$) is solved with zero initial iteration and ascending shift order. In \Cref{fig:ex3}, the left panel shows relative increment errors for the four methods at $n = 128$, while those related to warm restart and descending order for different $n$ are displayed in the right panel.

The results indicate that warm restart outperforms zero initial iteration by retaining information from smaller systems. For zero initial iteration, error curves begin at 1, and ascending shift order is preferred, consistent with the analysis in \Cref{sec:adi} and \Cref{subsec:order}, which prioritizes reducing dominant error components first. Conversely, with warm restart or a good initial guess, initial errors are smaller, making descending shift order more effective. Warm restart with descending order achieves errors below $\epsilon$ and stabilizes after approximately 30 iterations, whereas the ascending order variant fails to meet the stopping criterion until the final iteration, with incremental errors initially increasing. This indicates that high-frequency error terms, more prominent with warm restart, should be prioritized for reduction. These trends, depicted in \Cref{subfig:ex3relerr}, are consistent across all $n > 16$. Additionally, \Cref{subfig:ex3wr} reveals that initial errors decrease with increasing $n$ under warm restart, and the convergence rate with descending shift order remains stable throughout the solution process. For $n = 1024$, \eqref{optimalk} predicts over 80 shifts are required to reduce the error to $10^{-14}$. However, applying warm restart with descending shift order, fewer than 20 ADI iterations suffice to fully resolve the solution, reducing the computational cost by three-quarters.

\subsection{Factored ADI iterations}\label{subsec:fadi}
\begin{figure}[!t]
  \centering
  \subfloat{\label{subfig:ex4time}
  \includegraphics[width = 6.5cm]{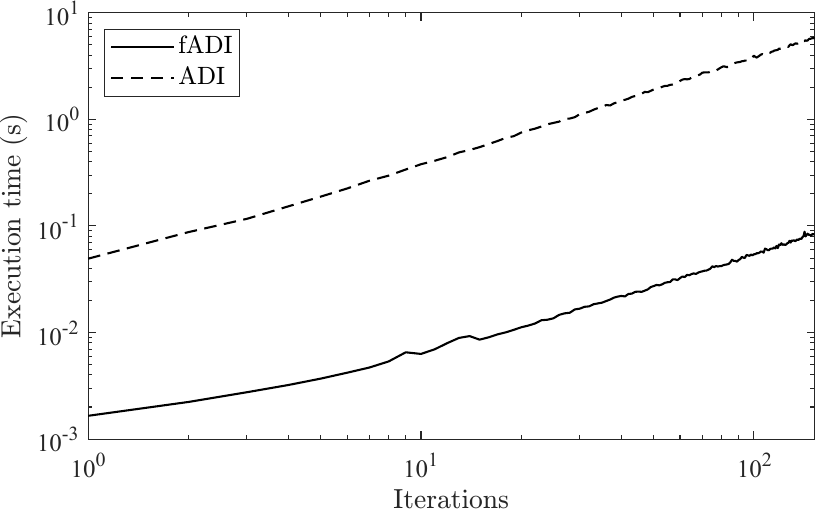}} \quad
  \subfloat{\label{subfig:ex4accuracy}
  \includegraphics[width = 6.5cm]{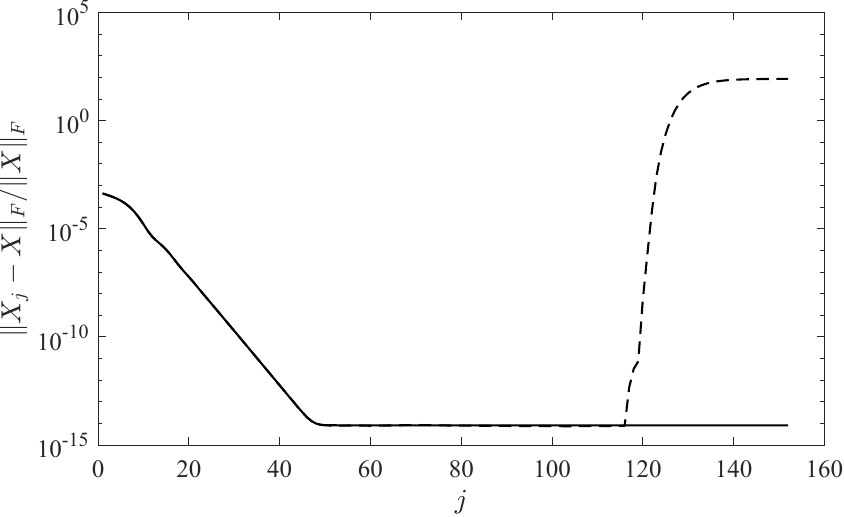}}
  \caption{Left: execution times against iterations for solving \eqref{fourth} using fADI \cref{fadiiter1} and ADI \cref{adiiter}. Right: relative errors of different methods.}
  \label{fig:ex4}
\end{figure}
The last example is a fourth order equation \eqref{fourth} with specified Dirichlet and Neumann boundary conditions so that the true solution is $u(x, y) = \sin(x^2)\exp(y^2)$. Note that $f = 16\exp(y^2)(x^4\sin(x^2) - 3x^2\cos(x^2) + 3y^2\sin(x^2) + y^4\sin(x^2))$ admits a low-rank function approximation. By subtracting a low-rank interpolation as described in \Cref{subsec:nonhomogeneous}, the term $F$ in \eqref{modelsyl} evidently possesses a low-rank factorization. For simplicity, we compute $F$ directly and apply singular value decomposition to obtain its low-rank form, revealing a numerical rank of 2 under a tolerance of $10^{-15}$. We compare the solving time and relative errors of the ADI method (\Cref{adiiter}) and the fADI method (\Cref{fadiiter2}), with results shown in \Cref{fig:ex4}. The truncation size is set to $n = 1000$, and 152 optimal shifts are used, determined via the analysis in \Cref{subsec:est}. The left panel of \Cref{fig:ex4} shows a speedup exceeding 70 times when using \eqref{fadiiter2}, highlighting the efficacy of the low-rank approach when $f$ is well-approximated by low-rank functions. In the right panel, the relative errors of both methods are indistinguishable up to 118 iterations, reflecting their inherent similarity. However, the ADI method exhibits a significant error spike near the 120th iteration, primarily due to shifts smaller than the machine epsilon of floating-point arithmetic. This issue is mitigated in practice, as the relative error stabilizes after 50 iterations, allowing early termination via the error-checking mechanism in \Cref{subsec:check}. In contrast, the fADI method remains robust against such tiny shifts, underscoring its stability for low-rank solutions.

\section{Concluding remarks}\label{sec:final}

We propose a method for solving the Poisson equation using the US method with optimal complexity, assuming sufficient smoothness in the underlying solution. This work provides valuable tools to fully leverage the speed and accuracy of the US method in two dimensions. We argue that the basis recombination introduced in \Cref{sec:us} offers a more suitable framework for extending the US method to higher dimensions and general equations. Preliminary research based on this approach has demonstrated promising results in both accuracy and complexity.

The optimality of the new method hinges on \Cref{thm:smoothness}, a factor often overlooked by other ADI solvers \citep{for, kno}. This theorem could also enhance those solvers, reducing their complexity to $\mathcal{O}(n^2)$. However, as these methods rely on Legendre polynomials, they are less competitive when solutions in Chebyshev polynomials are desired.

It is unsurprising that we focus solely on even-order problems, given their real spectra. Odd-order equations, with their complex spectra, pose a challenge for ADI methods. Nevertheless, such problems can be addressed by determining a region enclosing all eigenvalues \citep{wac2}. Initial experiments suggest that direct solvers, as presented here, may not outperform the Bartels--Stewart (B--S) method for these cases. Still, approximate solvers can be readily developed and serve as effective preconditioners when paired with iterative methods like GMRES for general equations. Further details will be explored in a forthcoming paper.

\section*{Acknowledgments}

We would like to acknowledge the heuristic discussion with Lu Cheng, Kuan Deng and Xiaolin Liu for \Cref{thm:smoothness}. We appreciate Kuan Xu for reading the manuscript and providing valuable feedback.

\bibliographystyle{abbrvnat}
\bibliography{optimalPoisson-bib}

\end{document}